\documentclass{amsart}
\usepackage{amscd,amssymb}
\usepackage{fullpage}
\usepackage{righttag}
\usepackage{pb-diagram,lamsarrow,pb-lams}

\DeclareMathOperator{\Aut}      {Aut}
\DeclareMathOperator{\cdepth}   {cdepth}
\DeclareMathOperator{\ord}      {ord}
\DeclareMathOperator{\Tot}      {Tot}
\DeclareMathOperator{\ann}      {ann}
\DeclareMathOperator{\Level}    {Level}
\DeclareMathOperator{\spf}      {spf}

\newcommand{\Fp}        {{\mathbb{F}_p}}
\newcommand{\Fpp}[1]    {{\mathbb{F}_p^{#1}}}
\newcommand{\Zpl}       {{\mathbb{Z}_{(p)}}}      
\newcommand{\Zp}        {{\mathbb{Z}_p}}          
\newcommand{\Zh}        {{\mathbb{Z}}}

\newcommand{\al}        {\alpha}
\newcommand{\bt}        {\beta} 
\newcommand{\tht}       {\theta}
\newcommand{\thth}      {\hat{\theta}}
\newcommand{\thb}       {\overline{\theta}}
\newcommand{\lm}        {\lambda}
\newcommand{\phih}      {\hat{\phi}}
\newcommand{\phb}       {\overline{\phi}}
\newcommand{\sg}        {\sigma}

\newcommand{\Sg}        {\Sigma}

\newcommand{\Ah}        {\widehat{A}}
\newcommand{\Eh}        {\widehat{E}}
\newcommand{\Fh}        {\widehat{F}}
\newcommand{\Ih}        {\widehat{I}}
\newcommand{\Jh}        {\widehat{J}}

\newcommand{\fps}[2]    {{#1 [\![ #2 ]\!]}} 
\newcommand{\tm}        {\times}
\newcommand{\ot}        {\otimes}
\newcommand{\Smash}     {\wedge}
\newcommand{\st}        {\;|\;}

\newcommand{\xra}       {\xrightarrow}
\newcommand{\xla}       {\xleftarrow}
\newcommand{\mra}       {\rightarrowtail}
\newcommand{\era}       {\twoheadrightarrow}

\newcommand{\ulm}       {\underline{\lambda}}
\newcommand{\ux}        {\underline{x}}
\newcommand{\BP}[1]     {{BP\langle #1\rangle}}
\newcommand{\Bb}        {\overline{B}}
\newcommand{\Cb}        {\overline{C}}
\newcommand{\cpi}       {{\mathbb{C}P^\infty}}
\newcommand{\invlim} {\operatornamewithlimits{\underset{\longleftarrow}{lim}}}
\newcommand{\OO}        {{\mathcal{O}}}
\newcommand{\GG}        {{\mathbb{G}}}
\newcommand{\ih}        {\hat{\imath}}
\newcommand{\mxi}       {\mathfrak{m}}
\newcommand{\catD}      {\mathcal{D}}
\newcommand{\catR}      {\mathcal{R}}
\newcommand{\Rb}        {\overline{R}}
\newcommand{\cp}        {{\mathbb{C}P}} 
\newcommand{\dps}       {\langle p\rangle(x)}
\newcommand{\xh}        {\hat{x}}

\renewcommand{\:}{\colon}

\newtheorem{theorem}{Theorem}[section]

\newtheorem{lemma}[theorem]{Lemma}
\newtheorem{proposition}[theorem]{Proposition}
\newtheorem{corollary}[theorem]{Corollary}
\theoremstyle{definition}

\newtheorem{definition}[theorem]{Definition}

\numberwithin{equation}{subsection}

\newenvironment{diag}{
 \renewcommand{\typeout}[1]{}
 \begin{displaymath}
 \begin{diagram}}{
 \end{diagram}
 \end{displaymath}} 

\begin{document}
\title{The {$\BP{n}$} cohomology of elementary Abelian groups}
\author{N.~P.~Strickland}
\bibliographystyle{abbrv}

\maketitle 

\section{Introduction}

In this paper we define three elements of a certain generalised
cohomology ring $\BP{m,n}^*BV_k$.  Here $m$, $n$ and $k$ are
nonnegative integers with $k+m\leq n+1$, there is a fixed prime $p$
not exhibited in the notation, and $V_k$ is an elementary Abelian
$p$-group of rank $k$.  We show that these elements are equal; this is
striking, because the three definitions are very different.  The
significance of our equation is not yet entirely clear, but it makes
contact with other work in the literature in a number of fascinating
ways.
\begin{enumerate}
 \item In the case where $n=1$, $m=0$ and $k=2$, our result is closely
  related to~\cite[Theorem 4.2]{anst:wpm}, which exhibits a relation
  in the connective complex $K$-theory group $kU^6B(\Zh/k)^2$ for all
  $k$.  That result is in turn a key part of a proof (at least in
  special cases) of the main result of~\cite{hoanst:esw}, which
  relates the Witten genus for spin manifolds~\cite{wi:ido} to the
  theory of elliptic curves and elliptic spectra.  The work described
  in the present paper started with an attempt to generalise this
  theorem in connective $K$-theory.
 \item Our methods also give an interesting filtration of the ring
  $\BP{n,m}^*BV_k$.  Again in the case $n=1$, $m=0$, $k=2$ our
  filtration is compatible with a splitting of $kU\Smash BV_2$
  discovered earlier by Ossa~\cite{os:cke}.  This was also a
  motivation for our investigations.  It would clearly be interesting
  to improve our filtration to some kind of stable splitting in the
  general case, but we have not succeeded as yet.
 \item One of our three definitions involves an iterated Bockstein map
  \[ q'_n := (\BP{0} \xra{q_1} \Sg^{2p-1} \BP{1} \xra{q_2} \ldots
     \xra{q_n} \Sg^{2(p^{n+1}-1)/(p-1)-n-2} \BP{n}).
  \]
  There is an analogous operation in Beilinson's motivic cohomology,
  which plays a central r\^{o}le in Voevodsky's proof of the Milnor
  conjecture in algebraic $K$-theory~\cite{vo:mc}.  We have no
  idea whether this fact is significant or not.
 \item By applying the $(n+2)$'th space functor to the map $q'_n$ of
  spectra, we get a map
  \[ r\:K(\Zpl,n+2) \xra{} \underline{\BP{n}}_{2(p^{n+1}-1)/(p-1)}. \]
  This map has appeared in a number of other places, for example the
  proof that suspension spectra are harmonic~\cite{hora:ssh}.  Moreover,
  Ravenel, Wilson and Yagita have shown that $BP^*r$ is surjective, and
  thus given a nice description of $BP^*K(\Zpl,n+2)$ \cite{rawiya:bpc}.
  Our results shed some interesting light on the nature of the map
  $r$.
 \item Our other two definitions make contact with the classical
  theory of Dickson invariants in $H^*(BV_k;\Fp)$ (see~\cite{wi:pdi}
  for an exposition).
 \item Our proofs use the theory of multiple level
  structures~\cite{grst:vlc}, which generalise Drinfel'd level
  structures on elliptic curves~\cite{dr:em}.  This theory was also
  central to the results of~\cite{grst:vlc}, which give a crude
  generalisation of the Hopkins-Kuhn-Ravenel generalised character
  theory for the rationalised Morava cohomology of classifying spaces
  of finite groups.
\end{enumerate}

There are other approaches to $\BP{n}_*BV_k$ and $\BP{n}^*BV_k$,
related to the Steinberg idempotents and the Conner-Floyd
conjecture~\cite{mi:pcf,jowi:bphi,jowiya:bphii}.  We have not yet
managed to find any fruitful interaction between our methods and
these, but it is clearly an interesting project to look for one.
 
\section{Statement of results}\label{sec-results}

We now state our results in more detail.  

Fix a prime $p$ and integers $m\le n$, and define $w=n+1-m$.  There is
a $BP$-algebra spectrum $E=\BP{m,n}$ with homotopy ring
$E^*=\Fp[v_m,\ldots,v_n]$ (or $E^*=\Zpl[v_1,\ldots,v_n]$ if $m=0$).
(We will recall a construction in Section~\ref{sec-spectra}.)  Note
that $w$ is the Krull dimension of $E^*$.  We use cohomological
gradings, so that $|v_k|=-2(p^k-1)<0$.  Write $V_k$ for the elementary
Abelian group $\Fpp{k}$.  We will give a filtration of $E^*BV_k$ by
ideals, such that each quotient is in a natural way a finitely
generated free module over a regular local ring.  We will show that
$E^*BV_k$ has no $v_m$-torsion when $k<w$.  When $k=w$, we show that
the $v_m$-torsion is annihilated by the ideal
$I_{n+1}=(v_m,\ldots,v_n)$, and that it is a free module on one
generator over the ring $\fps{\Fp}{x_0,\ldots,x_{w-1}}$.  We give
three very different formulae for this generator.

Recall that the \emph{depth} of a module $M^*$ over $R^*$ is the
largest $d$ such that there exists a sequence $\{a_0,\ldots,a_{d-1}\}$
in $R^*$ which is regular on $M^*$.  We will need an analogous but
different notion, which we now define.
\begin{definition}\label{defn-depth}
 If $M^*$ is a module over $E^*$, we define the \emph{chromatic depth}
 of $M^*$ to be the largest integer $d=\cdepth(M^*)$ such that the
 sequence $(v_m,v_{m+1},\ldots,v_{m+d-1})$ is regular on $M^*$.
 Clearly we have $0\le\cdepth(M^*)\le w$.
\end{definition}

We next introduce some abbreviated notation for rings which will
appear many times.
\begin{definition}\label{defn-rings}
 Given any algebra $R^*$ over $BP^*$, and integers $0\le i\le j\le k$,
 we write
 \begin{align*}
  P(R^*;j)    &= \fps{R^*}{x_t\st 0\le t<j}/([p](x_t)\st 0\le t<j)\\
  P(R^*;i,j)  &= \fps{R^*}{x_t\st i\le t<j}/([p](x_t)\st i\le t<j)\\
  P(R^*;i,j;k)&= \fps{R^*}{x_t\st i\le t<k}/([p](x_t)\st i\le t<j).
 \end{align*}
 Thus $P(R^*;j)=P(R^*;0,j)$ and $P(R^*;i,j)=P(R^*;i,j;j)$.
\end{definition}

The most basic fact is as follows.
\begin{proposition}\label{prop-basic}
 If $k\le w$ then we have 
 \[ E^*BV_k=P(E^*;k)=
     \fps{E^*}{x_0,\ldots,x_{k-1}}/([p](x_0),\ldots,[p](x_{k-1}))
 \]
 and this has chromatic depth $w-k$.
\end{proposition}
This is proved in Corollary~\ref{cor-depth-several} and
Corollary~\ref{cor-depth-geom}. 

In the extreme case $k=w$ we see that $E^*BV_w$ has chromatic depth
$0$, so it must have some $v_m$-torsion.  It is thus of interest to
see what the torsion subgroup is.  In order to state the answer, we
need some definitions.

\begin{definition}\label{defn-alpha}
 Given $\ulm\in\Fpp{j}$, we write
 $[\ulm](\ux)=[\lm_0](x_0)+_F\cdots +_F[\lm_{j-1}](x_{j-1})$.  We write
 \[ \phi_j(t) = \prod_{\ulm\in\Fpp{j}} (t-_F[\ulm](\ux))^{p^m}. \]
 We also define
 $\al=\al(m,n)=\prod_{j=0}^{w-1}\phi_j(x_j)\in E^*BV_w$.  Note that
 $\al$ is the product of the terms $[\ulm](\ux)^{p^m}$, where $\ulm$
 runs over elements of $\Fpp{w}$ whose last nonzero entry is one.
 The degree of $\al$ is $2p^m(p^w-1)/(p-1)$.
\end{definition}

In the following lemma, we regard $\BP{k}^*=\Zpl[v_1,\ldots,v_k]$ as a
subring of $\BP{n}^*$ in the obvious way.
\begin{lemma}\label{lem-pi}
 There are unique series $\pi_k(t)\in\fps{\BP{k}^*}{t}$ such that
 $[p](t)=\sum_{k=0}^nv_k\pi_k(t)\in\fps{\BP{n}^*}{t}$.  Moreover, the
 series $\pi_k(t)$ is divisible by $t^{p^k}$, and is equal to
 $t^{p^k}$ modulo $(p,t^{p^k+1})$ or modulo $I_{n+1}$.
\end{lemma}
\begin{proof}
 We give $t$ degree two, which makes everything homogeneous.  Write
 $\pi_0(t)=t$.  We know that $[p](t)=pt\pmod{t^2}$ so
 $[p](t)-v_0\pi_0(t)=\sum_{k>1}a_k t^k$ say.  Each $a_k$ has strictly
 negative cohomological degree (because $|t|=|[p](t)|=2$) so it lies
 in the ideal $(v_1,\ldots,v_n)$.  Let $f_k(t)$ be the sum of those
 monomials in $[p](t)-v_0\pi_0(t)$ that lie in $\fps{\BP{k}^*}{t}$,
 but do not lie in $\fps{\BP{k-1}^*}{t}$.  It is clear that $f_k(t)$
 is divisible by $v_k$, say $f_k(t)=v_k\pi_k(t)$.  This gives series
 $\pi_k(t)\in\fps{\BP{k}^*}{t}$ with $[p](t)=\sum_{k=0}^nv_k\pi_k(t)$,
 as required.  It is easy to see that they are unique.

 Note that the degree of $v_k\pi_k(t)$ is the same as the degree of
 $[p](t)$, which is $2$, so $\pi_k(t)$ has degree $2p^k$.  As $|t|=2$
 and $|v_i|<0$ for $i>0$ we see that $\pi_k(t)$ is divisible by
 $t^{p^k}$.  We thus have $\pi_k(t)=b_kt^{p^k}\pmod{t^{p^k+1}}$, say.
 One checks that $|b_k|=0$, so that $b_k\in\Zpl$.  Given that
 $[p](t)=v_kt^{p^k}\pmod{v_0,\ldots,v_{k-1},t^{p^k+1}}$, we see that
 $b_k=1\pmod{p}$, and thus $\pi_k(t)=t^{p^k}\pmod{p,t^{p^k+1}}$.  

 Similarly, the image of $\pi_k(t)$ modulo $I_{n+1}$ is a series in
 $\fps{\Fp}{t}$ of degree $2p^k$, so it has the form $c_kt^{p^k}$ for
 some $c_k\in\Fp$.  Because
 $[p](t)=v_kt^{p^k}\pmod{v_0,\ldots,v_{k-1},t^{p^k+1}}$, we see that
 $c_k=1$ and $\pi_k(t)=t^{p^k}\pmod{I_{n+1}}$.
\end{proof}

\begin{definition}\label{defn-alpha-prime}
 We define $\al'=\al'(m,n)$ to be the determinant of the square matrix
 with entries $\pi_i(x_j)$ for $m\le i\le n$ and $0\le j<w$.
\end{definition}

\begin{definition}\label{defn-alpha-sec}
 We shall see in Section~\ref{sec-spectra} that there are cofibrations
 of spectra
 \[ \Sg^{2p^n-2}\BP{m,n}\xra{v_n}\BP{m,n}\xra{}\BP{m,n-1}
     \xra{q_n}\Sg^{2p^n-1}\BP{m,n}.
 \]
 This also works for $n=m$ if we interpret $\BP{m,m-1}$ as the mod $p$
 Eilenberg-MacLane spectrum $H\Fp$.  We let $a_0,\ldots,a_{w-1}$ be
 the usual generators of $H^1(BV_w;\Fp)$ and we define
 \[ \al''=\al''(m,n)=q_n\cdots q_{m+1}q_m(a_0a_1\ldots a_{w-1}). \]
 It is not hard to see that when $g\in\Aut(V_w)$ we have
 $\det(g)\in\Fpp{\tm}$ and $g^*\al''=\det(g)\al''$ (which makes sense
 because $p\al''=0$).
\end{definition}

Our main result is as follows.
\begin{theorem}\label{thm-main}
 We have $\al=\al'=\pm\al''\in E^dBV_w$, where $d=2p^m(p^w-1)/(p-1)$.
 The annihilator of $v_m$ on $E^*BV_w$ is the same as the annihilator
 of $I_{n+1}$.  It is a free module of rank one over
 $E^*BV_w/I_{n+1}=P(\Fp;w)=\fps{\Fp}{x_0,\ldots,x_{w-1}}$ generated by
 $\al$.  It maps injectively to $H^*(BV_w;\Fp)$.
\end{theorem}
This will be proved after Proposition~\ref{prop-dickson}.  Note that
$\al'$ and $\al''$ appear to depend on the choice of generators $v_k$,
but in fact they do not, because they are equal to $\al$.  

The basic structure of the proof is as follows.  It is quite easy to
see that $I_{n+1}\al'=I_{n+1}\al''=0$ and that $\al$, $\al'$, and
$\al''$ all have the same degree.  We shall show in
Section~\ref{sec-spectra} that $q_k$ is compatible up to sign with the
Milnor Bockstein operation $Q_k$ in mod $p$ cohomology.  Given this,
classical arguments about Dickson invariants show that
$\al=\al'=\pm\al''\pmod{I_{n+1}}$.  We will show using an intricate
argument inspired by the theory of multiple level
structures~\cite{st:fsf,grst:vlc} that
$\ann(v_m)=\ann(I_{n+1})=(\al)$.  It follows for degree reasons that
$\al'=\lm\al$ for some $u\in\Fp$, and we find that $\lm=1$ by reducing
everything modulo $I_{n+1}$.  The same argument shows that
$\al''=\pm\al$.

\section{The spectra $\BP{m,n}$}\label{sec-spectra}

For brevity, we will write $MU$ for the $p$-local spectrum $MU_{(p)}$.
Recall that $MU$ can be made into a strictly commutative ring spectrum
(or ``commutative $S$-algebra'') in the foundational setting
of~\cite{ekmm:rma}, so we can construct a derived category
$\catD_{MU}$ of strict $MU$-modules, and a category $\catR_{MU}$ of
ring objects in $\catD_{MU}$ (referred to in~\cite{ekmm:rma} as
$MU$-ring spectra).  As usual, we let $BP^*$ be the largest quotient
ring of $MU^*$ over which the standard formal group law becomes
$p$-typical.  We know from~\cite{st:pmm} that (even when $p=2$) there
is a commutative ring $BP$ in $\catD_{MU}$ with homotopy ring $BP^*$,
and that this object is unique up to canonical isomorphism.

If $p>2$ we choose once and for all a sequence of elements 
$v_k\in BP_{2(p^k-1)}=BP^{-2(p^k-1)}$ such that
$[p](t)=v_kt^{p^k}\pmod{v_0,\ldots,v_{k-1},t^{p^k+1}}$.  In
particular, this forces $v_0=p$.  Two popular choices would be to take
$v_k$ to be the $k$'th Hazewinkel generator (so that
$[p](t)=\exp_F(pt)+_F\sum^F_{k>0}v_k t^{p^k}$) or the $k$'th Araki
generator (so that $[p](t)=\sum^F_{k\geq 0}v_kt^{p^k}$).  In any case,
we define $\BP{m,n}^*=BP^*/(v_i\st i<m\text{ or } i>n)$.  We know
from~\cite[Theorem 2.6]{st:pmm} that there is a commutative ring
$\BP{m,n}$ in $\catD_{MU}$ with homotopy ring $\BP{m,n}^*$, and
that this object is unique up to canonical isomorphism.

If $p=2$ we do not have so much choice about the sequence of $v$'s.
Nonetheless, we know from~\cite[Proposition 2.10]{st:pmm} that there
exist sequences of $v$'s for which the rings 
$\BP{n}^*=BP^*/(v_i\st i>n)$ can be realised as the homotopy ring of a
commutative ring object $\BP{n}$, which is unique up to canonical
isomorphism (here we are writing $\BP{n}$ for the spectrum called
$\BP{n}'$ in~\cite{st:pmm}).  Unfortunately, neither of the popular
choices listed above work in this context.  We fix such a sequence
once and for all.  By the proof of~\cite[Theorem 2.13]{st:pmm}, we
know that there is a central $\BP{n}$-algebra $\BP{m,n}$ with homotopy
ring $\BP{m,n}^*=\BP{n}^*/I_m$ and a derivation
$Q_{m-1}\:\BP{m,n}\xra{}\Sg^{2^m-1}\BP{m,n}$ such that
$ab-ba=v_mQ_{m-1}(a)Q_{m-1}(b)$.  It follows that when $X$ is a space
such that $\BP{m,n}^*X$ is concentrated in even degrees, the operation
$Q_{n-1}$ is trivial on the ring $\BP{n,m}^*X$, and thus this ring is
commutative.  If $\BP{m,n}$ and $BP'\langle m,n\rangle$ are two such
central $\BP{n}$-algebras then either there is a unique isomorphism
$\BP{m,n}\simeq BP'\langle m,n\rangle$, or there is a unique
isomorphism $\BP{m,n}^{\textup{op}}\simeq BP'\langle m,n\rangle$.

Next, recall from~\cite{st:pmm} that EKMM theory gives a map
$Q_k\:MU/v_k\xra{}\Sg^{2p^k-1}MU/v_k$ in $\catD_{MU}$ that is a
derivation for any product on $MU/v_k$.  By smashing this over $MU$
with $BP/(v_i\st i<m,i\neq k)$ we get a derivation
$Q_k\:P(m)\xra{}\Sg^{2p^k-1}P(m)$.  We can then smash this over $MU$
with $MU/(v_i\st i>n)$ and then with $MU/(v_m,\ldots,v_n)$ to get
compatible derivations on $\BP{m,n}$ and $H\Fp$.  It is also easy to
see that there is a canonical map
$q_n\:\BP{m,n-1}\xra{}\Sg^{2p^n-1}\BP{m,n}$ and a ring map
$\rho_n\:\BP{m,n}\xra{}\BP{m,n-1}$ fitting into a cofibre sequence 
\[ \Sg^{2p^n-2}\BP{m,n}\xra{v_n}\BP{m,n}\xra{\rho_n}
   \BP{m,n-1} \xra{q_n} \Sg^{2p^n-1}\BP{m,n}
\]
such that $\rho_n q_n=Q_n$.

We claim that our operation $Q_n$ on $H\Fp$ is the same as the Milnor
Bockstein operation $Q^M_n$ up to sign.  It is well-known that
derivations in mod $p$ cohomology are the same as primitive elements
in the Steenrod algebra $H\Fpp{*}H\Fp$, and that the space of
primitives in dimension $2p^n-1$ is spanned by $Q^M_n$, so we have
$Q_n=\lm Q^M_n$ for some $\lm\in\Fp$.  We need to show that $\lm=\pm
1$, so we need only compute one nontrivial instance of $Q_n$.  We
shall do this in the case $p>2$.  The case where $p=2$ and $n>0$
requires essentially only notational changes, and the case where $p=2$
and $n=0$ can be done with simple \emph{ad hoc} arguments.  Recall
that $H\Fpp{*}BV_1=\Fp[x]\ot E[a]$, where $a$ is the usual generator
of $H\Fpp{1}BV_1$ and $x=\bt a$, which is the image of the usual
generator of $MU^2\cpi$.  It is well-known that $Q^M_n(a)=x^{p^n}$.
Given this, our claim follows easily from
Proposition~\ref{prop-qktest} below.

\begin{proposition}\label{prop-qktest}
 Let $p$ be an odd prime, and let $C_p$ be the cyclic group of $p$'th
 roots of unity in $\mathbb{C}$, which we can identify with $V_1$.
 Write
 \[ X = \text{cofibre}(S^{2p^n-1}\xra{}S^{2p^n-1}/C_p), \]
 which is the $2p^n$-skeleton of $BC_p$.  Let $q_n$ be the Bockstein
 operation in the cofibration
 \[ \Sg^{2p^n-2} P(n) \xra{v_n} P(n) \xra{\rho}
     P(n+1) \xra{q_n} \Sg^{2p^n-1}P(n).
 \]
 Then there is a unique element $b\in P(n+1)^1X$ that hits the usual
 generator $a\in H\Fpp{1}X=H\Fpp{1}BC_p$.  Moreover,
 $q_n b=\pm x^{p^n}$, where $x$ is the image in $P(n)^2BC_p$ of the
 usual generator $x\in MU^2\cpi$.
\end{proposition}
The sign ambiguity could be resolved by a careful analysis of
conventions, which we do not have the patience to do.

\begin{proof}
 For brevity, write $m=p^n-1$ and $R=P(n)$ and $\Rb=P(n+1)$, so we
 have a cofibration
 \[ \Sg^{2m}R \xra{v_n} R \xra{\rho} \Rb \xra{q_n} \Sg^{2m+1}R.\]
 We also write
 \begin{align*}
   P   &= \cp^m = S^{2m+1}/S^1                                  \\
   L   &= \text{ tautological bundle over } P                   \\
   P^L &= \text{ Thom space of } L \simeq \cp^{m+1}             \\
   Y   &= S^{2m+1}/C_p                                          \\
   H   &= H\Fp.
 \end{align*}
 Note that $Y$ can also be thought of as the sphere bundle $S(L^p)$ in
 the $p$'th tensor power of $L$, or as the $(2m+1)$-skeleton of
 $BV_1=S^\infty/C_p$.  Note also that 
 \[ H^*X = P[x]\ot E[u]/(x^{m+2},a x^{m+1}) =
     \Fp\{1,a,x,\ldots,ax^m,x^{m+1}\}.
 \]

 An easy connectivity argument shows that there is a unique element
 $b\in\Rb^1X$ that hits $v\in H^1X$.  Indeed, let $F$ be the fibre of
 the map $\Rb\xra{}H$, so that $\pi_*F$ starts in dimension
 $2p^{n+1}-2$.  The bottom cell of $DX$ is in dimension $-2p^n$, so
 the bottom cell of $F\wedge DX$ is in dimension $d=2(p-1)p^n-2$.
 This is strictly greater than $1$ because $p>2$ and $n\ge 0$.  It
 follows easily that $\pi_1(DX\Smash\Rb)=\pi_1(DX\Smash H)$, or in
 other words $\Rb^1X=H^1X$, so there is a unique element $b$ as
 described.

 We next consider the diagram
 \begin{diag}
  \node[2]{P} \arrow{sw,l}{z_L} \arrow{se,t}{z_{L^p}} \\
  \node{P^L} \arrow[2]{e,b}{\phi} \node[2]{P^{L^p}}
 \end{diag}
 Here $z_L$ and $z_{L^p}$ are the zero-section inclusions, and $\phi$
 is obtained in the obvious way from the $p$'th power map of total
 spaces $E(L)\xra{}E(L^p)$.  There is also a corresponding diagram
 with $P$ replaced by $\cpi$.  By applying $R^*$ to this, we get a
 diagram
 \begin{diag}
  \node[2]{\fps{R^*}{x}} \\
  \node{\fps{R^*}{x}u_L}
  \arrow{ne,l}{z_L^*}
  \node[2]{\fps{R^*}{x}u_{L^p}}
  \arrow{nw,t}{z_{L^p}^*}
  \arrow[2]{w,b}{\phi^*}
 \end{diag}
 Here $u_L$ and $u_{L^p}$ are the Thom classes of $L$ and $L^p$; the
 corresponding Euler classes are $x$ and $[p](x)$.  It is well-known
 that $z_L^*(f(x)u_L)=f(x)x$ and $z_{L^p}^*(f(x)u_{L^p})=f(x)[p](x)$,
 and it follows easily (because $\dps=[p](x)/x$ is not a zero-divisor
 in $\fps{R^*}{x}$) that
 \[ \phi^*(f(x)u_{L^p})=f(x)\dps u_L. \]
 By mapping the $\cp^m$ diagram into the $\cp^\infty$ diagram, we see
 that the same fomula holds for the map 
 \[ \phi^*\:\tilde{R}^*P^{L^p}=(R^*[x]/x^{m+1}).u_{L^p} \xra{}
            (R^*[x]/x^{m+1}).u_L = \tilde{R}^*P^L.
 \]
 In this context, however, we simply have $\dps=v_nx^m$, so
 \[ \phi^*u_{L^p} = v_n x^m u_L. \]

 We next apply the octahedral axiom to the maps
 $S^{2m+1}\xra{}Y\xra{}P$ (which comes down to replacing the maps by
 cofibrations and using the fact that
 $P/Y=(P/S^{2m+1})/(Y/S^{2m+1})$).  As $S^{2m+1}=S(L)$ and
 $Y=S^{2m+1}/C_p=S(L^p)$, we see that the cofibres of
 $S^{2m+1}\xra{}P$ and $Y\xra{}P$ are homeomorphic to $P^L$ and
 $P^{L^p}$.  With these identifications, the map $S^{2m+1}\xra{}Y$ is
 the $p$'th power map $S(L)\xra{}S(L^p)$, and thus the induced map
 $P^L\xra{}P^{L^p}$ is just $\phi$.  We therefore have an octahedral
 diagram
 \[
 \setlength{\unitlength}{0.14em}
 \begin{picture}(140,110)(0,-10)
  \put( 18,11){\vector( 2, 3){18}}
  \put( 70,11){\vector( 2, 3){18}}
  \put( 62,77){\vector(-2,-3){18}}
  \put( 44,38){\vector( 2,-3){18}}
  \put( 96,38){\vector( 2,-3){18}}
  \put( 88,50){\vector(-2, 3){18}}
  \put(112, 5){\vector(-1, 0){40}}
  \put( 60, 5){\vector(-1, 0){40}}
  \put( 52,44){\vector( 1, 0){34}}
  \put( 10,11){\line(-2, 3){ 8}}
  \put(  2,23){\line( 2, 3){40}}
  \put( 42,83){\vector( 1, 0){18}}
  \put( 73,83){\line( 1, 0){17}}
  \put( 90,83){\line( 2,-3){40}}
  \put(130,23){\vector(-2,-3){ 8}}
  \put( 27,25){\circle{3}}
  \put( 53,64){\circle{3}}
  \put( 92, 5){\circle{3}}
  \put( 14, 5){\makebox(0,0){$X$}}
  \put( 66, 5){\makebox(0,0){$Y$}}
  \put(118, 5){\makebox(0,0){$P^{L^p}$}}
  \put( 40,44){\makebox(0,0){$S^{2m+1}$}}
  \put( 92,44){\makebox(0,0){$P$}}
  \put( 66,83){\makebox(0,0){$P^L$}}
  \put( 88,23){\makebox(0,0){$\pi_{L^p}$}}
  \put( 46,25){\makebox(0,0){$\phi$}}
  \put( 66,48){\makebox(0,0){$\pi_L$}}
  \put( 86,64){\makebox(0,0){$z_L$}}
  \put( 46,64){\makebox(0,0){$\tau_L$}}
  \put( 20,25){\makebox(0,0){$\tau_{L^p}$}}
  \put( 37, 0){\makebox(0,0){$j$}}
  \put( 95, 0){\makebox(0,0){$e$}}
  \put(112,25){\makebox(0,0){$z_{L^p}$}}
  \put( 15,59){\makebox(0,0){$r$}}
  \put(117,59){\makebox(0,0){$\phi$}}
 \end{picture}
 \]
 (A circled arrow
 $U\longrightarrow\hspace{-1.3em}\circ\hspace{0.8em}V$ means a map
 $U\xra{}\Sigma V$.  The diagram can be made to look more like an
 octahedron by lifting up the outer three vertices and drawing in an
 extra arrow to represent the composite
 $je\:P^{L^p}\longrightarrow\hspace{-1.3em}\circ\hspace{0.8em}X$.)  In
 particular, we see that the stable fibre of $\phi$ is just $X$.

 We next claim that the map
 $e^*\:\widetilde{H}^1Y\xra{}\widetilde{H}^2P^{L^p}$ is an
 isomorphism, and sends $a$ to the Thom class $u_{L^p}$.  To see this,
 let $M$ be the restriction of $L$ to the basepoint in $P$, so $M$ is
 just a one-dimensional complex vector space.  The bottom cell of
 $Y=S(L^p)$ is just the circle $S(M^p)$, and the bottom cell of
 $P^{L^p}$ is the one-point compactification of $M^p$, written
 $S^{M^p}$.  The restriction of $e$ to this bottom cell is the
 standard homeomorphism of $S^{M^p}$ with the unreduced suspension of
 $S(M^p)$, followed by the projection to the reduced suspension.  The
 claim follows easily from this.

 We now consider the following diagram.
 \begin{diag}
  \node{\Sg^{-2}P^L}
  \arrow{s,l}{x^m u_L}
  \arrow{e,t}{\phi}
  \node{\Sg^{-2}P^{L^p}}
  \arrow{s,l}{u_{L^p}}
  \arrow{e,t}{je}
  \node{\Sg^{-1}X}
  \arrow{e,t}{r}
  \arrow{s,r,..}{c}
  \node{\Sg^{-1}P^L}
  \arrow{s,r}{x^m u_L} \\
  \node{\Sg^{2m}R}
  \arrow{e,b}{v_n}
  \node{R}
  \arrow{e,b}{\rho}
  \node{\Rb}
  \arrow{e,b}{q_n}
  \node{\Sg^{2m+1}R}
 \end{diag}
 We see from the octahedron that the top line is a cofibration, and
 the bottom line is a cofibration by construction.  The left hand
 square commutes because $\phi^*u_{L^p}=v_nx^mu_L$.  It follows that
 there exists a map $c\:\Sg^{-1}X\xra{}\Rb$ (i.e. $c\in \Rb^1X$)
 making the whole diagram commute.

 From our discussion of $e^*$ in cohomology and the commutativity of
 the middle square, we see that the image of $c$ in $H^1X$ is just
 $a$.  By the uniqueness of $b$, we deduce that $c=b$.  Thus, the
 commutativity of the right hand square tells us that $q_n b$ is the
 image of $x^mu_L\in R^{2m+2}P^L$ in $R^*X$.  Under the usual
 identification of $P^L$ with $\cp^{m+1}$, the element $x^mu_L$
 becomes $x^{m+1}$, and the map $X\xra{}\cp^{m+1}$ is a restriction of
 the usual map $BC_p\xra{}\cpi$.  It follows that $q_n b=x^{m+1}$ as
 claimed.
\end{proof}

\section{Commutative algebra}\label{sec-comm-alg}

In this section, we recall some basic ideas from commutative algebra.

We will need to be a little more careful than is usual about the
relationship between graded and ungraded rings.  For us, a graded ring
$R^*$ will mean a sequence of Abelian groups $R^k$ (for $k\in\Zh$)
with product maps $R^i\ot R^j\xra{}R^{i+j}$ with the usual properties.
We will assume that $R^k=0$ when $k$ is odd, and that the product is
commutative.  This will apply to all rings that we consider except for
$H^*(BV_k;\Fp)$, and in that case we will not need to use the results
of this section.  It is common to identify a graded ring $R^*$ with
the ungraded ring $\bigoplus_k R^k$.  We shall not do this, for the
following reason.  There is an obvious way to interpret the expression
$R^*=\fps{\Zpl[v_1,\ldots,v_n]}{x_0,\ldots,x_{k-1}}$ as a graded ring
(with $|v_k|=-2(p^k-1)$ and $|x_i|=2$).  Explicitly, $R^k$ is the set
of expressions $\sum_\al a_\al x^\al$, where
$\al=(\al_0,\ldots,\al_{k-1})$ is a multiindex, $|\al|=\sum_i\al_i$
and $a_\al\in\BP{n}^{-2|\al|}$.  With this interpretation, the ring
$\bigoplus_kR^k$ is not the same as the ungraded ring
$R'=\fps{\Zpl[v_1,\ldots,v_n]}{x_0,\ldots,x_{k-1}}$ (it does not
contain $\sum_{i\ge 0}x_0^i$, for example).  The set $\prod_kR_k$ is
different again.  It is not clear how many of the nice ring-theoretic
properties of $R'$ are shared by $\bigoplus_kR^k$.  For this reason,
we prefer to work with graded rings as defined above.  Most theorems
for ungraded rings have graded counterparts, which are proved by a
straightforward adaptation of the ungraded proofs.  We will outline
the results that we need, leaving most of the task of adaptation to
the reader.

We will only consider homogeneous elements of $R^*$, in other words
elements of $R^k$ for various $k$.  The word ``module'' will always
mean ``graded module'', and similarly for ideals.  

Given an ideal $I^*$ in $R^*$, let $(I^N)^k$ denote the part of the
$N$'th power of $I$ in degree $k$.  The cosets $a+(I^N)^k$ (for $a\in
R^k$) form a basis for a topology on $R^k$.  We say that a topology of
this form is a \emph{linear topology} on $R^*$, and we say that $I^*$
is an \emph{ideal of definition}.  If $I^*$ is one ideal of
definition, it is clear that another ideal $J^*$ is an ideal of
definition for the same topology if $I^*$ contains a power of $J^*$
and \emph{vice versa}.  We say that a linear topology is
\emph{complete} if $R^k=\invlim_N R^k/(I^N)^k$ for all $k$.  We say
that a homogeneous element is \emph{topologically nilpotent} if some
power of it lies in $I^*$.  In the Noetherian case, it is not hard to
see that the topologically nilpotent elements form an ideal of
definition. 

If $R^*$ is a quotient of $\fps{E^*}{x_0,\ldots,x_{k-1}}$ then we give
$R^*$ the complete linear topology defined by the ideal
$(x_0,\ldots,x_{k-1})$.  In this context, we can consider the ungraded
ring $\Tot(R^*)=\invlim_N\bigoplus_k(R/I^N)^k$ as a substitute for
$\bigoplus_kR^k$.  This avoids the difficulty mentioned above: if
\[ R^*=\fps{\Fp[v_m,\ldots,v_n]}{x_0,\ldots,x_{k-1}} \] 
(in the usual graded sense) then
\[ \Tot(R^*)=\fps{\Fp[v_m,\ldots,v_n]}{x_0,\ldots,x_{k-1}} \]
(in the usual ungraded sense).  However, the relationship between
properties of $R^*$ and those of $\Tot(R^*)$ is not as close as one
might like, so we will not stress this point of view.

We will need the following version of the Weierstrass preparation
theorem. 
\begin{proposition}\label{prop-weierstrass}
 Let $R^*$ be as above.  Let $y=\sum_{k\ge 0}a_kx^k\in\fps{R^*}{x}$ be
 a homogeneous element such that $a_i$ is topologically nilpotent for
 $i<d$, and $a_d$ is a unit.  Then $\fps{R^*}{x}$ is freely generated
 by $\{x^i\st i<d\}$ as a module over $\fps{R^*}{y}$, and thus
 $\fps{R^*}{x}/y$ is freely generated by $\{x^i\st i<d\}$ as a module
 over $R^*$.
\end{proposition}
\begin{proof}
 We may assume that $a_d=1$.  For any $m\ge 0$, we can write $m=ld+k$
 with $l\ge 0$ and $0\le k<d$, and then write $w_m=x^ky^l$.
 Let $I^*$ be the ideal of topological nilpotents, so that
 $w_m=x^m\pmod{I^*,x^{m+1}}$.  For any
 $R^*$-module $M^*$ that is complete with respect to $I^*$, we can
 define a map
 \[ \tht_M\:\prod_{m\ge 0} M^* \xra{} \fps{M^*}{x}
 \]
 by $\tht_M(b)=\sum_mb_mw_m$.  Suppose that $I^*M^*=0$, and that
 $c=\sum_{m\ge 0}c_mx^m\in\fps{M^*}{x}$.  Given the form of $w_m$ mod
 $I^*$, we see easily by induction that there is a unique sequence of
 $b_m$'s such that $c=\sum_{k=0}^mb_mw_m\pmod{x^{m+1}}$ and thus that
 $\tht_M$ is an isomorphism.  Moreover, if we have a pair of modules
 $N^*\leq M^*$ such that $\tht_N$ and $\tht_{M/N}$ are isomorphisms,
 then so is $\tht_M$ (by a five-lemma argument).  It follows by
 induction that $\tht_{M/I^k}$ is iso for all $k$, and thus by taking
 inverse limits that $\tht_M$ is iso for all $M^*$.  In particular,
 $\tht_R$ is an isomorphism.  This means that for any series
 $c(x)\in\fps{R^*}{x}$, there are unique series
 $b_0(y),\ldots,b_{d-1}(y)\in\fps{R^*}{y}$ such that $c(x)=\sum_i
 b_i(y)x^i$, which proves the proposition.
\end{proof}

We say that a graded ring $R^*$ is \emph{local} if it has only one
maximal ideal, or equivalently if it has an ideal such that every
homogeneous element in the complement is invertible.  It is clear that
$E^*$ is local in this sense, with maximal ideal $(v_m,\ldots,v_n)$.
Similarly, any quotient of the ring $\fps{E^*}{x_0,\ldots,x_{k-1}}$ is
local, with maximal ideal $(v_m,\ldots,v_n,x_0,\ldots,x_{k-1})$.

We say that a graded ring $R^*$ is Noetherian if every ideal is
generated by a finite set of homogeneous elements.  Simple adaptations
of the usual arguments in an ungraded 
context~\cite[Theorem 3.3]{ma:crt} show that any quotient of
$\fps{E^*}{x_0,\ldots,x_{k-1}}$ is Noetherian.

We say that a graded ring $R^*$ is a \emph{domain} if the product of
two nonzero homogeneous elements of $R^*$ is nonzero.  This will hold
if $\Tot(R^*)$ is an ungraded domain, but unfortunately we have not
been able to prove the converse.  We say that an ideal $P^*$ in $R^*$
is \emph{prime} if $R^*/P^*$ is a domain, and define the \emph{Krull
  dimension} of $R^*$ to be the largest integer $d$ such that there
exists a chain $P^*_0<\ldots<P^*_d$ of prime ideals in $R^*$.

We say that a Noetherian graded local ring $R^*$ of Krull dimension
$d$ is \emph{regular} if there is a sequence of $d$ homogeneous
elements that generates the maximal ideal.  Such a sequence is called
a \emph{regular system of parameters}; it is necessarily a regular
sequence.  Simple adaptations of the usual arguments in an ungraded
context~\cite[Theorem 15.4]{ma:crt} show that
$R^*=\fps{E^*}{x_0,\ldots,x_{k-1}}$ has dimension $w+k$, so the
sequence $\{v_m,\ldots,v_n,x_0,\ldots,x_{k-1}\}$ is a regular sequence
of parameters and $R^*$ is regular.  More generally, if $R^*$ is any
graded regular local ring and $x$ is a homogeneous indeterminate then
$\fps{R^*}{x}$ is again a graded regular local ring.

By graded versions of~\cite[Theorems 14.3 and 20.3]{ma:crt}, we see
that a graded regular local ring is an integral domain, with unique
factorisation for homogeneous elements.

The following result is a graded version of~\cite[Theorem
14.2]{ma:crt}, and can be proved in the same way.
\begin{theorem}\label{thm-reg-sys}
 Let $R^*$ be a graded regular local ring of dimension $d$, and
 $\{x_0,\ldots,x_{k-1}\}$ a sequence of homogeneous elements of
 $R^*$.  Then the following are equivalent:
 \begin{itemize}
 \item[(a)] $\{x_0,\ldots,x_{k-1}\}$ is a subset of a regular system
  of parameters.
 \item[(b)] The images of $\{x_0,\ldots,x_{k-1}\}$ in $\mxi/\mxi^2$
  (where $\mxi$ is the unique maximal ideal of $R^*$) are linearly
  independent over the graded field $R^*/\mxi$.
 \item[(c)] $R^*/(x_0,\ldots,x_{k-1})$ is a graded regular local ring
  of dimension $d-k$.
 \end{itemize}
\end{theorem}

We can also show that two weaker conditions are equivalent to each
other. 
\begin{theorem}\label{thm-reg-seq}
 Let $R^*$ be a graded regular local ring of dimension $d$, and
 $\{x_0,\ldots,x_{k-1}\}$ a sequence of homogeneous elements of
 $R^*$.  Then the following are equivalent:
 \begin{itemize}
 \item[(a)] $\{x_0,\ldots,x_{k-1}\}$ is a regular sequence.
 \item[(b)] $R^*/(x_0,\ldots,x_{k-1})$ has Krull dimension $d-k$.
 \end{itemize}
\end{theorem}
\begin{proof}
 This is a graded version of the equivalence (1)$\Leftrightarrow$(3)
 in part~(iii) of~\cite[Theorem 17.4]{ma:crt}, taking into account the
 equation $\text{ht}(I^*)=\dim(R^*)-\dim(R^*/I^*)$ from part~(i) of
 that theorem, and the fact that regular local rings are
 Cohen-Macaulay.  
\end{proof}

\section{The structure of $E^*BV_k$}

\begin{proposition}
 Let $R^*$ be an algebra over $E^*$ of chromatic depth $d>0$.  Then
 $[p](x)$ is not a zero-divisor in $\fps{R^*}{x}$, and the ring
 $\fps{R^*}{x}/[p](x)$ has chromatic depth at least $d-1$.
\end{proposition}

\begin{proof}
 Suppose that $0\neq f(x)\in\fps{R^*}{x}$, say
 $f(x)=ax^k\pmod{x^{k+1}}$ with $0\neq a\in R^*$.  As $R^*$ has
 nonzero chromatic depth, we see that $v_ma\neq 0$ and
 $[p](x)f(x)=v_max^{k+p^m}\pmod{x^{1+k+p^m}}$ so $[p](x)f(x)\neq 0$.
 This shows that $[p](x)$ is not a zero-divisor in $\fps{R^*}{x}$.

 Now suppose that $d>1$, so that $\{v_m,v_{m+1}\}$ is regular on
 $R^*$.  Consider a series $f(x)\in\fps{R^*}{x}$ such that $v_mf(x)$
 is divisible by $[p](x)$.  We claim that $f(x)$ is divisible by
 $[p](x)$.  Indeed, suppose that $v_mf(x)=[p](x)g(x)$.  We then have
 $[p](x)g(x)=0$ in $\fps{(R^*/v_m)}{x}$.  However, $v_{m+1}$ is not a
 zero divisor in $R^*/v_m$, so $[p](x)$ is not a zero divisor in
 $\fps{(R^*/v_m)}{x}$ by the previous paragraph.  Thus $g(x)=0$ in
 $\fps{R^*/v_m}{x}$, or $g(x)=v_mh(x)$ in $\fps{R^*}{x}$ say.  This
 means that $v_m(f(x)-[p](x)h(x))=0$ in $\fps{R^*}{x}$, but $v_m$ is
 not a zero-divisor so $f(x)=[p](x)h(x)$ as claimed.  The conclusion
 is that $v_m$ is regular on $\fps{R^*}{x}/[p](x)$.  If $d>2$ then we
 can replace $R^*$ by $R^*/v_m$ and $m$ by $m+1$ and $d$ by $d-1$ and
 run the same argument again.  It follows by induction that
 $\fps{R^*}{x}/[p](x)$ has chromatic depth at least $d-1$.
\end{proof}

\begin{corollary}\label{cor-depth-several}
 Let $R^*$ be an algebra over $E^*$ of chromatic depth $d>0$.  Then
 for $j\le d$ and $k\ge j$ the sequence $[p](x_0),\ldots,[p](x_{j-1})$
 is regular in $\fps{R^*}{x_0,\ldots,x_{k-1}}$ and the quotient
 $P(R^*;j;k)$ has chromatic depth at least $d-j$.
\end{corollary}
\begin{proof}
 This follows easily from the proposition, using the obvious fact that
 $\fps{S^*}{x}$ has the same chromatic depth as $S^*$ for any algebra
 $S^*$ over $\BP{m,n}^*$.
\end{proof}

\begin{corollary}\label{cor-depth-geom}
 For $k\le w$ we have $E^*BV_k=P(E^*;k)$, and this has chromatic depth
 at least $w-k$.
\end{corollary}
\begin{proof}
 Write $B=B\Zh/p$ and $Z=\cpi$, so $E^*Z=\fps{E^*}{x}$
 and $B$ is a circle bundle over $Z$ with Chern class $[p](x)$.  Thus
 $B^i\tm Z^{k-i}$ is a circle bundle over $B^{i-1}\tm Z^{k+1-i}$ with
 Chern class $[p](x_i)$.  This gives a long exact Gysin sequence 
 \[  \cdots \xla{} E^*(B^i\tm Z^{k-i}) \xla{}
     E^*(B^{i-1}\tm Z^{k+1-i}) \xla{[p](x_i)}
     E^{*-2}(B^{i-1}\tm Z^{k+1-i}) \xla{} \cdots.
 \]
 It follows easily by induction from Corollary~\ref{cor-depth-several}
 that these sequences are short exact, and that
 $E^*(B^i\tm Z^{k-i})=P(E^*;i;k)$.  The case $i=k$ gives the
 corollary. 
\end{proof}

\begin{definition}
 If $k\le w$, we define $A(k)^*=A(m,n,k)^*$ to be the largest quotient
 ring of $R^*=P(E^*;k)$ over which the series
 $\phi_k(t)=\prod_{\ulm\in\Fpp{k}}(t-_F[\ulm](\ux))^{p^m}$ divides
 $[p](t)$.  We also write $\psi_k(t)=[p](t)/\phi_k(t)$.

 In more detail, we note that $t-_F[\ulm](\ux)$ is a unit multiple of
 $t-[\ulm](\ux)$, so $\phi_k(t)$ is a unit multiple of a
 monic polynomial of degree $p^{m+k}$, whose lower coefficients are
 topologically nilpotent.  It follows from
 Proposition~\ref{prop-weierstrass} that $\fps{R^*}{t}/\phi_k(t)$ is a
 free module over $R^*$ on generators $1,t,\ldots,t^{p^{m+k}-1}$.  In
 particular, we can write
 $[p](t)=\sum_{i=0}^{p^{m+k}-1}c_it^i\pmod{[p](t)}$, for uniquely
 defined coefficients $c_i\in R^*$.  We define $A(k)^*$ to be the
 quotient ring $R^*/(c_0,\ldots,c_{p^{m+k}-1})$.  It also follows from
 Proposition~\ref{prop-weierstrass} that $\phi_k(t)$ is not a
 zero-divisor in $\fps{A(k)^*}{t}$, so there is a unique series
 $\psi_k(t)\in \fps{A(k)^*}{t}$ such that $[p](t)=\phi_k(t)\psi_k(t)$.
\end{definition}

\begin{theorem}\label{thm-Ak}
 For $k\le w$ we have
 \begin{enumerate}
 \item[(a)] $A(k)^*$ is a regular local ring in the graded sense.
 \item[(b)] There is a unique formal group law $F_k$ over $A(k)^*$ for
  which $\phi_k(t)$ is a homomorphism $F\xra{}F_k$.  
 \item[(c)] There is a unique series $\tht_k(t)=t\thb_k(t)$ over
  $A(k)^*$ such that $[p]_F(t)=\tht_k(\phi_k(t))$.
 \item[(d)] If $k>0$, we have
  $A(k)^*=\fps{A(k-1)^*}{x_{k-1}}/\psi_{k-1}(x_{k-1})$.
 \end{enumerate}
 Moreover, when $k=w$ we have
 $A(w)^*=P(\Fp;w)=\fps{\Fp}{x_0,\ldots,x_{w-1}}$ and $I_{n+1}=0$ in
 $A(w)^*$.
\end{theorem}
\begin{proof}
 Suppose that $k<w$ and that (a)$\ldots$(d) hold up to stage $k$
 (which is trivial for $k=0$).  Define
 $B^*=\fps{A(k)^*}{x_k}/\psi_k(x_k)$.

 We first claim that $\psi_k(x_k)=0$ in $A(k+1)^*$.  Indeed, we
 have 
 \[ \phi_{k+1}(t)=\phi_k(t)\prod_{i=1}^{p-1}\phi_k(t-_F[i](x_k)), \]
 and this divides $[p](t)=\phi_k(t)\psi_k(t)$ over $A(k+1)^*$.  As
 $\phi_k(t)$ is a unit multiple of a monic polynomial, it is not a
 zero divisor in $\fps{A(k+1)^*}{t}$, so we conclude that
 $\prod_{i=1}^{p-1}\phi_k(t-_F[i](x_k))$ divides $\psi_k(t)$ in
 $\fps{A(k+1)^*}{t}$.  We now set $t=x_k$ to conclude that
 $\psi_k(x_k)=0$ in $A(k+1)^*$ as claimed.

 It is also clear that $\phi_k(t)$ divides $\phi_{k+1}(t)$, so it
 divides $[p](t)$ over $A(k+1)^*$, so $A(k+1)^*$ is an algebra over
 $A(k)^*$.  This means that $A(k+1)^*$ is a quotient of $B^*$ in a
 natural way.

 Next, we claim that $B^*$ is a regular local ring in the graded
 sense.  For this, we write $C^*=\fps{E^*}{x_0,\ldots,x_k}$.  This is
 clearly a regular local ring, and
 $B^*=C^*/(\psi_0(x_0),\ldots,\psi_k(x_k))$.  By
 Theorem~\ref{thm-reg-sys} it suffices to check that the list
 \[ L=\{\psi_0(x_0),\ldots,\psi_k(x_k),
      x_0,\ldots,x_k,v_{m+k+1},\ldots,v_n\}
 \]
 is a regular system of parameters for $C^*$.  The length of $L$ is
 $n+k-m+2$ which is the same as the Krull dimension of $C^*$, so it is
 enough to check that $L$ generates the maximal ideal of $C^*$, or
 equivalently that $\Cb^*=C^*/(L)=\Fp$, or equivalently that
 $v_m,\ldots,v_{m+k}$ vanish in $\Cb^*$, or equivalently that $[p](t)$
 is divisible by $t^{p^{m+k+1}}$ over $\Cb^*$.  However, we know that
 $C^*$ is an algebra over $A(k)^*$ so $[p](t)=\phi_k(t)\psi_k(t)$ over
 $\Cb^*$.  Moreover, as $x_0,\ldots,x_{k-1}$ vanish in $\Cb^*$ we see
 from the definition that $\phi_k(t)=t^{p^{k+m}}$ in $\fps{\Cb^*}{t}$.
 This means that $[p](t)=t^{p^{k+m}}\psi_k(t)$ over $\Cb^*$, so
 $v_0=\ldots=v_{k+m-1}=0$ in $\Cb^*$ and $v_{k+m}=\psi_k(0)$.  On the
 other hand, we have $x_k=0$ and $\psi_k(x_k)=0$ so $\psi_k(0)=0$ so
 $v_{k+m}=0$ as required.  This completes the proof that $B^*$ is a
 regular local ring.

 It is clear from the above that
 $\{x_0,\ldots,x_k,v_{m+k+1},\ldots,v_n\}$ is a regular system of
 parameters for $B^*$, and thus that the Krull dimension of $B^*$ is
 $w$.  As $B^*$ is an integral domain, it follows that any proper
 quotient of $B^*$ has dimension strictly less than $w$.  Thus, if we
 can show that $A(k+1)^*$ has dimension $w$, we can deduce that the
 quotient map $B^*\era A(k+1)^*$ is an isomorphism.

 For the rest of the argument we need to separate the cases $k+1=w$
 and $k+1<w$.  If $k+1=w$ we argue as follows.  It follows from the
 definition of $A(w)^*$ that $A(w)^*/I_{n+1}$ is the largest quotient
 of $\fps{E^*}{x_0,\ldots,x_{w-1}}/I_{n+1}=P(\Fp;w)$ over which
 $\phi_w(t)$ divides $[p](t)$.  However, in this context $[p](t)=0$
 which is automatically divisible by $\phi_w(t)$, so
 $A(w)^*/I_{n+1}=P(\Fp;w)$.  This has the same Krull dimension as
 $B^*$, so we must have $B^*=A(w)^*=A(w)^*/I_{n+1}$.  In particular,
 we see that $I_{n+1}=0$ in $A(w)^*$, and thus that the formal group
 law $F$ becomes $F(s,t)=s+t$ in $\fps{A(w)^*}{s,t}$.  It follows that
 parts~(b) and~(c) of the theorem hold with $F_w=F$ and $\tht_w=0$.

 We now consider the case $k+1<w$, so that $k<w-1$.  For this, we use
 the theory of multiple level structures developed in~\cite[Section
 5]{grst:vlc}; we will assume that the reader is familiar with this.

 Let $\Eh^*$ be the graded ring
 $(E^*[u,u^{-1}]/(u^{p^n-1}-v_n))^\wedge_{I_n}$.  Write
 $u_k=v_ku^{1-p^k}\in\Eh^0$.  It is easy to see that
 $\Eh^0=\fps{\Fp}{u_m,\ldots,u_{n-1}}$, or
 $\fps{\Zp}{u_1,\ldots,u_{n-1}}$ if $m=0$.  Moreover, we have
 $\Eh^*=\Eh^0[u,u^{-1}]$. 

 Let $F$ be the usual formal group law over $E^*$, and define
 $\Fh(s,t)=uF(s/u,t/u)$, so that $\Fh$ is a formal group law over
 $\Eh^0$.  Let $\GG$ be the associated formal group over
 $X=\spf(\Eh^0)$, which has height $n$ and strict height $m$.  This
 puts us in the context studied in~\cite{grst:vlc}.

 Now write $J=I_n+(x_0,\ldots,x_k)\leq A(k+1)^*$, and
 $\Ah(k+1)^*=(A(k+1)^*[u,u^{-1}]/(u^{p^n-1}-v_n))^\wedge_J$.  It is
 clear from the definitions that $\Ah(k+1)^0$ is just the ring
 $\OO_{\Level_m(V_{k+1}^*,\GG)}$ which classifies $p^m$-fold
 level-$V_{k+1}^*$ structures on $\GG$, and
 $\Ah(k+1)^*=\Ah(k+1)^0[u,u^{-1}]$.  We know from~\cite[Theorem
 5.6]{grst:vlc} that $\Ah(k+1)^0$ is an integral domain and a finitely
 generated free module over $\Eh^0$.  It is nonzero because
 $k+1\le n-m$.  It follows that $\Ah(k+1)^*$ is a graded domain and a
 finitely generated free module on homogeneous generators over
 $\Eh^*$. 

 We now have a diagram as follows.
 \begin{diag}
  \node{E^*} \arrow{e,t}{i} \arrow{s,l}{c}
  \node{A(k+1)^*} \arrow{s,r}{d} \\
  \node{\Eh^*} \arrow{e,b}{\ih} \node{\Ah(k+1)^*.} 
 \end{diag}
 For $m\le t\le n$ we let $\Ih_t$ be the ideal in $\Eh^*$ generated by
 $\{v_m,\ldots,v_{t-1}\}$.  This is clearly prime, with
 $\Ih_t\cap{}E^*=I_t$ and $0=I_m<\ldots<I_n$.  As $\Ah(k+1)^*$ is
 integral over $\Eh^*$, the going-up theorem~\cite[Theorem
 9.3]{ma:crt} gives us a chain $0=\Jh_m<\ldots<\Jh_n$ of primes in
 $\Ah(k+1)^*$ such that $\Jh_t\cap\Eh^*=\Ih_t$.  Let $J_t$ be the
 preimage of $\Jh_t$ under the map $d\:A(k+1)^*\xra{}\Ah(k+1)^*$.  It
 is clear that the preimage of this under the map
 $i\:E^*\xra{}A(k+1)^*$ is just $I_t$, so we have a chain of strict
 inclusions of prime ideals $0=J_m<\ldots<J_n$.  Let $J_{n+1}$ be the
 maximal ideal in $A(k+1)^*$.  As $v_n$ is a unit in $\Ah(k+1)^*$ it
 is clear that $v_n\not\in J_n$, but $v_n\in J_{n+1}$ so we have a
 chain of strict inclusions $0=J_m<\ldots<J_{n+1}$.  This shows that
 $A(k+1)^*$ has Krull dimension $w$, so $B^*=A(k+1)^*$ as explained
 earlier.
 
 We now define $\xh_i=ux_i\in\Ah(k+1)^0$, and for $\ulm\in\Fpp{k+1}$ we
 define
 \[ [\ulm](\xh) = 
    [\lm_0]_{\Fh}(\xh_0)+_{\Fh}\cdots+_{\Fh}[\lm_k]_{\Fh}(\xh_k).
 \]
 We also write 
 \[ \phih_{k+1}(t)=
     \prod_{\ulm}(t-_{\Fh}[\ulm](\xh))^{p^m}=
     u^{p^{m+k}}\phi_{k+1}(t/u)\in\fps{\Ah(k+1)^0}{t}.
 \]
 Clearly, the map $\ulm\mapsto [\ulm](\ux)$ defines a $p^m$-fold
 level-$V_{k+1}^*$ structure on $\GG$ over $\spf(\Ah(k+1)^0)$.  If we
 compose this with the projection $V_m^*\tm V_{k+1}^*\xra{}V_{k+1}^*$
 we get an ordinary ($1$-fold) level-$(V_m\tm V_{k+1})^*$ structure.
 It follows from~\cite[Proposition 32 and Corollary 33]{st:fsf} that
 $\phih_{k+1}(t)$ is a coordinate on a quotient formal group of $\GG$,
 and that the kernel of the quotient map is contained in the kernel of
 $p\:\GG\xra{}\GG$.  This means that there is a unique formal group
 law $\Fh_{k+1}$ defined over $\Ah(k+1)^0$ such that $\phih_{k+1}$ is
 a homomorphism $\Fh\xra{}\Fh_{k+1}$, and that there is a unique power
 series $\thth_{k+1}(t)$ defined over $\Ah(k+1)^0$ such that
 $[p]_{\Fh}(t)=\thth_{k+1}(\phih_{k+1}(t))$.  We define
 $F_{k+1}(s,t)=u^{-p^k}\Fh_{k+1}(u^{p^k}s,u^{p^k}t)$ and
 $\tht_k(t)=u^{-1}\thth_{k+1}(u^{p^k}t)$.  It is easy to deduce that
 $F_{k+1}$ is the unique formal group law over $\Ah(k+1)^*$ such that
 $\phi_{k+1}$ is a homomorphism $F\xra{}F_{k+1}$, and that
 $\tht_{k+1}$ is the unique series such that
 $[p]_F(t)=\tht_{k+1}(\phi_{k+1}(t))$.

 We next remark that $A(k+1)^*$ is a Noetherian domain, so it is
 not hard to see that the maps 
 \[ A(k+1)^*\xra{}v_n^{-1}A(k+1)^* \xra{}
     (v_n^{-1}A(k+1)^*)^\wedge_J = \Ah(k+1)^* 
 \]
 are injective.

 Now write $s'=\phi_{k+1}(s)$ and $t'=\phi_{k+1}(t)$.  We know from
 Proposition~\ref{prop-weierstrass} that $\fps{A(k+1)^*}{s,t}$ is a
 free module over the subring $\fps{A(k+1)^*}{s',t'}$ on generators
 $s^it^j$ for $0\le i,j<p^{k+1}$.  We can thus write
 $\phi_{k+1}(F(s,t))=\sum_{i,j}F_{i,j}(s',t')s^it^j$ for uniquely
 determined series $F_{i,j}$.  Similarly, $\fps{\Ah(k+1)^*}{s,t}$ is a
 free module over $\fps{\Ah(k+1)^*}{s',t'}$ on $\{s^it^j\}$, so the
 equation $\phi_{k+1}(F(s,t))=\sum_{i,j}F_{i,j}(s',t')s^it^j$ is the
 \emph{unique} way to write $\phi_{k+1}(F(s,t))$ in terms of the
 generators $s^it^j$.  On the other hand, we also have
 $\phi_{k+1}(F(s,t))=F_{k+1}(s',t')\in\fps{\Ah(k+1)^*}{s,t}$.  It
 follows that $F_{i,j}=0$ for $(i,j)\neq(0,0)$ and that $F_k=F_{0,0}$,
 so the series $F_k$ is actually defined over $A(k+1)^*$ rather than
 $\Ah(k+1)^*$.  It is clearly the unique formal group law over
 $A(k+1)^*$ for which the series $\phi_{k+1}$ is a homomorphism
 $F\xra{}F_k$.  Similarly, we see that $\tht_{k+1}(t)$ is defined over
 $A(k+1)^*$, and it is the unique series over $A(k+1)^*$ for which
 $[p](t)=\tht_{k+1}(\phi_{k+1}(t))$.  By putting $t=0$ we see that
 $\tht_{k+1}(t)$ is divisible by $t$, so it can be written as
 $t\thb_{k+1}(t)$.  This completes our induction step.
\end{proof}

\begin{proposition}\label{prop-depth-Ak}
 The ring $A(k)^*$ has chromatic depth at least $w-k$.
\end{proposition}
\begin{proof}
 The claim is that the sequence $\{v_m,\ldots,v_{n-k}\}$ (of length
 $w-k$) is regular on $A(k)^*$.  By Theorem~\ref{thm-reg-seq}, it is
 enough to check that the quotient $A(k)^*/(v_m,\ldots,v_{n-k})$ has
 dimension at most $\dim(A(k)^*)-(w-k)=k$.  Let $B^*$ be a graded
 integral domain which is a quotient of $A(k)^*/(v_m,\ldots,v_{n-k})$.
 It is enough to show that every such $B^*$ has dimension at most $k$,
 and thus enough to show that the maximal ideal $\mxi$ of $B^*$ needs
 at most $k$ generators.  From now on we work in $B^*$.

 Write $W=\{\ulm\in\Fpp{k}\st [\ulm](\ux)=0\}$.  This is a subgroup of
 $V_k^*$, of dimension $d$ say.  Note that $\Aut(\Fpp{k})=GL_k(\Fp)$
 acts on $A(k)^*$ in a natural way.  After applying a suitable element
 of this group, we may assume that $W$ is the evident copy of
 $\Fpp{d}$ in $\Fpp{k}$, spanned by the first $d$ standard basis
 vectors, so that $x_0=\ldots=x_{d-1}=0$ in $B$.  We write $U$ for the
 space spanned by the remaining standard basis vectors, so that
 $\Fpp{k}=W\oplus U$.  Define
 $\phb(t)=\prod_{\ulm\in U}(t-_F[\ulm](\ux))$, so that
 $\phi_k(t)=\phb(t)^{p^{m+d}}$.  As $B^*$ is an integral domain, we
 see that
 \[ \phb'(0)=\prod_{\ulm\in U\setminus 0}[\ulm](\ux)\neq 0 \]
 and thus that $\ord_t\phi_k(t)=p^{m+d}$, where $\ord_tf(t)$ means the
 largest integer $N$ such that $t^N$ divides $f(t)$.

 As $v_0=\ldots=v_{n-k}=0$ in $B^*$ we see that
 $\ord_t[p](t)\geq p^{n+1-k}$.  On the other hand, we have
 $\ord_t[p](t)=\ord_t\tht_k(\phb(t)^{p^{m+d}})=p^{m+d}\ord_s\tht_k(s)$.   
 Thus $\ord_s\tht_k(s)\ge p^{w-d-k}$.

 Now consider the list
 \[ L = \{x_d,\ldots,x_{k-1},v_{n+1-d},\ldots,v_n\}, \]
 so that $L$ has length $k$.  It will be enough to show that
 $\Bb=B/(L)=\Fp$, or equivalently that $x_0=\ldots=x_k=0$ in
 $\Bb$ and $v_0=\ldots=v_n=0$ in $\Bb$.  We already have
 $x_0=\ldots=x_{d-1}=0$ in $B$ and the remaining $x$'s are in $L$ so
 all $x$'s vanish in $\Bb$, as required.  We next note that
 $\phi_k(t)$ becomes $t^{p^{m+k}}$ over $\Bb$, and
 $[p](t)=\tht_k(\phi_k(t))$, and $\ord_s(\tht_k(s))\ge p^{w-d-k}$ so
 $[p](t)$ has height at least $m+k+w-d-k=n+1-d$ over $\Bb$.  This
 means that $v_0=\ldots=v_{n-d}=0$ over $\Bb$ and the remaining
 $v$'s are in $L$ so they vanish in $\Bb$ also.  Thus $(L)=\mxi$,
 and $\mxi$ needs only $k$ generators, as required.
\end{proof}

We leave the proof of the following simple lemmas to the reader.
\begin{lemma}\label{lem-quot}
 Let $R^*$ be a ring, and let $\phi,\psi$ be elements of $R^*$ such
 that $\phi$ is not a zero-divisor.  Then the annihilator of $\phi$ in
 $R^*/\phi\psi$ is generated by $\psi$, and thus the ideal generated
 by $\phi$ in $R^*/\phi\psi$ is a free module of rank one over
 $R^*/\psi$.  \qed
\end{lemma}
\begin{lemma}
 For any elements $\phi$ and $\chi$ in any ring $R^*$, the annihilator
 $\ann(\phi\chi,R^*)$ in $R^*$ is the preimage of
 $\ann(\phi,R^*/\ann(\chi))$ under the quotient map
 $R^*\xra{}R^*/\ann(\chi)$.  \qed
\end{lemma}

\begin{definition}
 We write $\phi_i$ for $\phi_i(x_i)$ and 
 $\chi_j=\prod_{i<j}\phi_i\in E^*BV_j\leq E^*BV_w$.
 We also write $\psi_i$ for $\psi_i(x_i)=\thb_i(\phi_i(x_i))$, so that
 $\phi_i\psi_i=[p](x_i)$.  We also have
 $E^*BV_k/(\psi_i\st i<j)=P(A(j)^*;j,k)$ by Theorem~\ref{thm-Ak}(d).
\end{definition}

\begin{proposition}\label{prop-ann-psi}
 If $j\leq k\leq w$ then the annihilator of $\chi_j$ in $E^*BV_k$ is
 precisely the ideal $(\psi_0,\ldots,\psi_{j-1})$.  Thus, the ideal
 generated by $\chi_j$ is a free module of rank one over
 $P(A(j)^*;j,k)$.
\end{proposition}
\begin{proof}
 The proposition is trivial when $j=0$ (so $\chi_j=1$).  We may thus
 assume the statement for $\chi_j$ and prove the one for $\chi_{j+1}$.
 As $\chi_{j+1}=\phi_j\chi_j$, we see that the annihilator of
 $\chi_{j+1}$ is the preimage of the annihilator of $\phi_j$ in
 $S^*=E^*BV_k/\ann(\chi_j)=P(A(j)^*;j,k)$.  We write
 $R^*=P(A(j)^*;j+1,k)$, so that $S^*=\fps{R^*}{x_j}/\phi_j\psi_j$.  It
 follows from Lemma~\ref{lem-quot} that the annihilator of $\phi_j$ in
 $S^*$ is just the ideal generated by $\psi_j$.  Thus, we have
 \begin{align*} 
  \ann(\chi_{j+1},E^*BV_k) &= \ann(\chi_j,E^*BV_k) + (\psi_j)   \\
   &= (\psi_0,\ldots,\psi_{j-1}) + (\psi_j)                     \\
   &= (\psi_0,\ldots,\psi_j),
 \end{align*}
 as required.
\end{proof}

The following corollary is immediate.
\begin{corollary}\label{cor-quotients}
 If $0\leq j<k\leq w$ then the quotient $(\chi_j)/(\chi_{j+1})$ is
 isomorphic as a module over $E^*BV_k$ to
 $\fps{P(A(j)^*;j+1,k)}{x_j}/\phi_j$, which is (by
 Proposition~\ref{prop-weierstrass}) a free module of rank $p^{mj}$
 over $P(A(j)^*;j+1,k)$. \qed
\end{corollary}

\begin{proposition}\label{prop-main}
 The ideal of $v_m$-torsion elements in $E^*BV_w$ is a free module
 over the ring $E^*BV_w/I_{n+1}=\fps{\Fp}{x_0,\ldots,x_{w-1}}$,
 generated by $\al=\chi_w$.  The map
 $\ann(v_m)\mra E^*BV_w\era E^*BV_w/I_{n+1}$ is injective.
\end{proposition}
\begin{proof}
 Take $k=w$ in Corollary~\ref{cor-quotients}.  We know from
 Proposition~\ref{prop-depth-Ak} that $A(j)^*$ has chromatic depth at
 least $w-j$, and it follows from Corollary~\ref{cor-depth-several}
 that $P(A(j)^*;j+1,w)$ has chromatic depth at least
 $(w-j)-(w-(j+1))=1$, so $v_m$ is regular on $(\chi_j)/(\chi_{j+1})$.
 It follows easily that $v_m$ is regular on $E^*BV_w/(\chi_w)$, and
 thus that the $v_m$-torsion is contained in the ideal $(\chi_w)$.  We
 know from Proposition~\ref{prop-ann-psi} that the annihilator of
 $\chi_w$ is the same as the kernel of the map $E^*BV_w\era A(w)^*$.
 We also know from Theorem~\ref{thm-Ak} that
 $A(w)^*=E^*BV_w/I_{n+1}=\fps{\Fp}{x_0,\ldots,x_{w-1}}$.  It follows
 immediately that the $v_m$-torsion is precisely the ideal generated
 by $\chi_w$ and that it is a free module over $E^*BV_w/I_{n+1}$.  It
 is easy to see from the definitions that $\chi_w=\al$.  Moreover,
 modulo $I_{n+1}$ we see that $\chi_w$ is just the product of all
 linear polynomials of the form $\lm_0x_0+\ldots+\lm_{j-1}x_{j-1}+x_j$
 for some $j<w$.  In particular, it is nonzero.  As $E^*BV_w/I_{n+1}$
 is an integral domain, it follows easily that the map
 $(E^*BV_w/I_{n+1}).\chi_w\mra E^*BV_w\era E^*BV_w/I_{n+1}$ is
 injective.
\end{proof}

Our proof of Theorem~\ref{thm-main} relies on some facts from the
classical theory of Dickson invariants~\cite{wi:pdi}.  For ease of
reference, we give a swift proof of what we need.
\begin{proposition}\label{prop-dickson}
 In $H^*(BV_k;\Fp)=E[a_0,\ldots,a_{k-1}]\ot\Fp[x_0,\ldots,x_{k-1}]$ we
 define
 \[ \bt = \prod_u \sum_i \lm_i x_i, \]
 where the product runs over nonzero sequences
 $(\lm_0,\ldots,\lm_{k-1})\in\Fpp{k}$ whose last nonzero entry is one.
 We also define
 \begin{align*}
  \bt'     &= \det_{i,j}(x_i^{p^j})                       \\
  \bt''_m  &= Q_{m+k-1}\ldots Q_m(a_0a_1\ldots a_{k-1}) 
 \end{align*}
 Then $\bt''_m=(\bt')^{p^m}=\bt^{p^m}$.
\end{proposition}
\begin{proof}
 Using column operations and the fact that
 $(u+v)^{p^j}=u^{p^j}+v^{p^j}\pmod{p}$, we see easily that
 $g^*\bt'=\det(g)\bt'$ for all $g\in\Aut(V_k)$.  It is also clear that
 the element $a_0a_1\ldots a_{k-1}$ transforms in the same way, and
 thus that $g^*\bt''_m=\det(g)\bt''_m$ for all $g$ and $m$.  It is
 immediate from the definition that $\bt'$ lies in
 $\Fp[x_0,\ldots,x_{k-1}]$ and that it is divisible by $x_0$.  It is
 well-known that $Q_i$ is a derivation with $Q_i(a_j)=x_j^{p^i}$ and
 $Q_i(x_j)=0$, and from this it follows that $\bt''_0$ is a sum of
 terms of the form $\pm\prod_{i<k}x_i^{p^{\sg(i)}}$ for various
 permutations $\sg$.  In particular, we see that $\bt''_0$ lies in
 $\Fp[x_0,\ldots,x_{k-1}]$ and that it is divisible by $x_0$.  As
 $\bt'$ and $\bt''_0$ are divisible by $x_0$ and invariant under
 $SL_k(\Fp)$, we see that they are both divisible by each of the terms
 $\sum_i\lm_ix_i$ in the product formula for $\bt$.  As these terms are
 inequivalent irreducibles and $\Fp[x_0,\ldots,x_{k-1}]$ has unique
 factorisation, we see that $\bt''_0$ and $\bt'$ are divisible by
 $\bt$.  It is easy to check that $\bt''_0/\bt$ and $\bt'/\bt$ have
 degree zero, so they lie in $\Fp$.  By comparing coefficients for the
 monomial $\prod_{i<k}x_i^{p^i}$, we see that $\bt=\bt'=\bt''_0$.

 Finally, we can define an (ungraded) ring map
 $F\:H^*(BV_k;\Fp)\xra{}H^*(BV_k;\Fp)$ by $F(a_i)=a_i$ and
 $F(x_i)=x_i^p$.  We then have $F\circ Q_i=Q_{i+1}\circ F$, from which
 it follows easily that $(\bt''_0)^{p^m}=F^m\bt''_0=\bt''_m$, which
 completes the proof of the proposition.
\end{proof}

\begin{proof}[Proof of Theorem~\ref{thm-main}]
 Given Proposition~\ref{prop-main}, all that is left is to prove that
 $\al=\al'=\pm\al''$.  

 We first show that $v_i\al'=v_i\al''=0$ for $i=m,\ldots,n$.

 Recall that $\al'$ is the determinant of the matrix with entries
 $\pi_i(x_j)$ for $m\le i\le n$ and $0\le j<w$.  As we work modulo
 $I_m$ we have $[p](t)=\sum_{i=m}^n v_i\pi_i(t)$, so
 $v_i\pi_i(x_j)=-\sum_{k\neq i}^nv_k\pi_k(x_j)$ for all $j$.  This shows
 that when we multiply the $i$'th row of our matrix by $v_i$, it
 becomes a linear combination of the other rows, so the determinant
 becomes zero.  It follows that $v_i$ annihilates the determinant of
 the original matrix, in other words $v_i\al'=0$.

 Next, choose elements $v_i\in\pi_*MU$ lifting our elements
 $v_i\in\pi_*BP$.  We then have natural maps
 $v_i\:\Sg^{2p^i-2}M\xra{}M$ for all $M\in\catD_{MU}$.  Working in
 $\catD_{MU}$ we see that $v_iq_i=0$ and thus
 $v_iq_n\ldots q_{m+1}q_m=0$.  After applying the forgetful functor to
 the category of spectra, we conclude that $v_i\al''=0$.

 In view of Proposition~\ref{prop-main}, it is now enough to check
 that $\al$, $\al'$ and $\pm\al''$ have the same image modulo
 $I_{n+1}$, or equivalently the same image in $H^*(BV_w;\Fp)$.  One
 can check from the definitions that these are the same as the
 elements $\bt^{p^m}$, $(\bt')^{p^m}$ and $\pm\bt''_m$ of
 Proposition~\ref{prop-dickson}, so they are all the same, as
 required.
\end{proof}

We conclude with a finer filtration of $E^*BV_w$.
\begin{proposition}\label{prop-filtration}
 The ring $E^*BV_w$ admits a finite filtration by ideals, whose
 quotients are finitely generated free modules over regular local
 rings of dimension $w$.
\end{proposition}
This could be made more explicit but the bookkeeping would be
tedious.  
\begin{proof}
 We actually prove that the ring $P(A(i)^*;j,k)$ admits such a
 filtration whenever $i\leq j\leq k\leq w$.  The case when $i=j=0$ and
 $k=w$ gives the proposition.  When $j=k$ we have
 $P(A(i)^*;j,k)=A(i)^*$ so the claim follows from
 Theorem~\ref{thm-Ak}, so we work by induction on $k-j$.  We have
 $P(A(i)^*;j,k)\simeq P(A(i)^*;i,k+i-j)$, so we may assume that $j=i$.
 Note that
 $P(A(i)^*;i,k)=\fps{P(A(i)^*;i+1,k)}{x_i}/\phi_i(x_i)\psi_i(x_i)$, so
 we can apply Lemma~\ref{lem-quot}.  This gives us a two stage
 filtration of $P(A(i)^*;i,k)$ in which one quotient is
 $P(A(i+1)^*;i+1,k)$ and the other is a finitely generated free module
 over $P(A(i)^*;i+1,k)$.  The proposition follows by induction.
\end{proof}

In the case where $n=1$ and $m=0$, the reduced $E$-cohomology of
$BV_1$ is a free module of rank one over $A(1)^*=\fps{E^*}{x}/\dps$,
and we have $A(2)^*=\fps{\Fp}{x_0,x_1}=H\Fpp{*}(\cpi\tm\cpi)$.  Ossa
has shown that $E\Smash BV_{2+}$ splits as an $E$-module as a wedge of
copies of $E$, $E\Smash BV_1$ and $H\Fp\Smash(\cpi\tm\cpi)_+$.  (He
actually works with the connective $K$-theory spectrum $kU$, but it is
well-known that this splits $p$-locally as a wedge of copies of $E$.
He also works with $BV_1\Smash BV_1$ rather than $BV_{2+}$ but again
the translation is trivial.)  One can check that the induced splitting
of $E^*BV_2$ splits our filtration of $E^*BV_2$.

In more general cases, it is unclear what happens.  The most plausible
idea seems to be that there should be an $\BP{m,n}$-algebra spectrum
$A(m,n,k)$ with homotopy ring $A(m,n,k)^*$ and that 
$\BP{m,n}\Smash BV_w$ should be a finite wedge of spectra of the form
$A(m',n,k')$ for various $m'$ and $k'$.  However, much work remains to
be done in this direction.


\end{document}